\documentclass[12pt,reqno,a4paper]{amsart}
\title[Sign-graded Posets]{Sign-graded posets, unimodality of 
$W$-polynomials and the Charney-Davis Conjecture}

\thanks{Part of this work was  
financed by the EC's IHRP Programme, within the Research 
Training Network 
``Algebraic Combinatorics in Europe'', grant HPRN-CT-2001-00272, while the 
author was at Universit\'a di Roma ``Tor Vergata'', Rome, Italy.}

\usepackage[english]{babel}
\usepackage{amsmath}
\usepackage{amsthm,amssymb,latexsym}
\usepackage{t1enc}
\usepackage{epic}
\usepackage{eepic}
\usepackage{xypic}
\usepackage{mathrsfs}

\numberwithin{equation}{section}
\newtheorem{proposition}{Proposition}[section]
\newtheorem{lemma}[proposition]{Lemma}
\newtheorem{corollary}[proposition]{Corollary}
\newtheorem{theorem}[proposition]{Theorem}
\newtheorem{conjecture}[proposition]{Conjecture}
\theoremstyle{definition}
\newtheorem{definition}[proposition]{Definition}

\newcommand{\N}{\mathbb{N}}
\newcommand{\om}{\omega}

\newcommand{\Z}{\mathbb{Z}}
\newcommand{\J}{\mathcal{L}}
\newcommand{\R}{\mathbb{R}}

\newcommand{\Q}{\mathbb{Q}}

\newcommand{\Pa}{\mathcal{A}}

\newcommand{\e}{\ar@{-}}
\newcommand{\de}{\ar@{.}}
\newcommand{\n}{*=0{\bullet}}

\def\ch#1,#2,{\binom#1#2}
\def\emm#1,{{\em #1}}
\def\ba#1,{\overline{#1}}

\def\gen#1,{\langle #1 \rangle}
\def\geno#1,{\langle #1 \rangle_{\infty}}

\def\newop#1{\expandafter\def\csname #1\endcsname{\mathop{\rm
#1}\nolimits}}

\newop{span}
\newop{mode}
\newop{min}
\newop{max}
\newop{deg}
\newop{des}
\newop{Con}
\newop{Ann}
\newop{dim}
\newop{Int}
\newop{Im}
\newop{ker}
\newop{mult}
\newop{ex}

\usepackage{graphicx}

\begin{document}
\maketitle
\begin{center}
{\sc Petter Br\"and\'en}\\
\medskip
{Department of Mathematics,}\\
{Chalmers University of Technology and G\"oteborg University}\\
{S-412~96  G\"oteborg, Sweden}\\
{\tt branden@math.chalmers.se}
\end{center}
\bigskip

\begin{abstract}
We generalize the notion of graded posets to what we call sign-graded 
(labeled) posets. We prove that the $W$-polynomial of a sign-graded 
poset is symmetric and unimodal. This extends a recent result 
of Reiner and Welker who proved it for graded posets by associating 
a simplicial polytopal sphere to 
each graded poset. By proving that the $W$-polynomials of sign-graded 
posets has the right sign at $-1$, we are able to prove the Charney-Davis 
Conjecture for these spheres (whenever they are flag).
\end{abstract} 
\nocite{*}\bibliographystyle{plain}

\thispagestyle{empty}
\section{Introduction and preliminaries}
Recently Reiner and Welker \cite{reinerwelker} proved that the 
$W$-polynomial of 
a graded poset (partially ordered set) $P$ has 
unimodal coefficients. They proved this by associating to $P$ a simplicial 
polytopal sphere, $\Delta_{eq}(P)$, whose $h$-polynomial is the 
$W$-polynomial of $P$, 
and invoking the $g$-theorem for simplicial polytopes 
(see \cite{stanleyg,stanleyalg}). Whenever 
this sphere is flag, i.e., its minimal non-faces all have cardinality two, 
they noted that the Neggers-Stanley Conjecture implies 
the Charney-Davis Conjecture for $\Delta_{eq}(P)$. In this 
paper we give a different proof of the unimodality 
of $W$-polynomials of graded posets, and we also prove the Charney-Davis 
Conjecture for $\Delta_{eq}(P)$ (whenever it is flag). 
We prove it by studying a family 
of labeled posets, which we call sign-graded posets, of which the class of 
graded naturally labeled posets is a sub-class.     
 
In this paper all posets will be finite and non-empty. For undefined 
terminology on posets 
we refer the reader to \cite{stanley1}. We denote the cardinality of 
a poset $P$ with the letter $p$. Let $P$ be a poset   
and let 
$\omega : P \rightarrow \{1,2,\ldots ,p\}$ be a bijection. The pair 
$(P,\omega)$ is called a {\em labeled poset}. If $\omega$ is order-preserving 
then $(P,\omega)$ is said to be {\em naturally labeled}.  
A $(P,\om)$-{\em partition} is a map $\sigma : P \rightarrow \{1,2,3,\ldots\}$ 
such that
\begin{itemize}
\item $\sigma$ is order reversing, that is, if $x \leq y$ then 
      $\sigma(x) \geq \sigma(y)$, 
\item if $x < y$ and $\om(x) > \om(y)$ then $\sigma(x) > \sigma(y)$.
\end{itemize}
The theory of $(P,\om)$-partitions was developed by Stanley in 
\cite{stanleythesis}.  
The number of $(P,\om)$-partitions $\sigma$ with largest part 
at most $n$   
is a polynomial of degree $p$ in $n$ called the {\em order polynomial} of 
$(P,\om)$ and is denoted $\Omega(P,\om;n)$.   
The $W$-polynomial of $(P,\om)$ is defined by 
$$
\sum_{n \geq 0}\Omega(P,\om;n+1)t^n = \frac{W(P,\om;t)}{(1-t)^{p+1}}. 
$$
The set, $\J(P,\om)$, 
of permutations $\om(x_1), \om(x_2), \ldots, \om(x_{p})$ where 
$x_1, x_2, \ldots, x_{p}$ is a linear extension of $P$ is called 
the {\em Jordan-H\"older set} of $(P,\om)$.  
A {\em descent} in a permutation $\pi= \pi_1\pi_2\cdots \pi_p$ is 
an index $1 \leq i \leq p-1$ such that $\pi_i > \pi_{i+1}$. The 
number of descents in $\pi$ is denoted $\des(\pi)$.  
A fundamental result on $(P,\om)$-partitions, see \cite{stanleythesis},    
is that the $W$-polynomial can 
be written as 
$$
W(P,\om;t)=\sum_{\pi \in \J(P,\om)}t^{\des(\pi)}. 
$$
The Neggers-Stanley Conjecture is the following:
\begin{conjecture}[Neggers-Stanley]
For any labeled poset $(P,\om)$ the polynomial $W(P,\om;t)$ has 
only real zeros.
\end{conjecture}
This was first conjectured by Neggers \cite{neggers} in 1978 for 
natural labelings and by Stanley in 1986 for arbitrary labelings. The 
conjecture has been proved for some special cases, see 
\cite{branden1, brentithesis, reinerwelker, wagnernc} for the state of the 
art.  
If a polynomial has only real non-positive zeros then its coefficients 
form a unimodal sequence. For the $W$-polynomials of graded 
posets unimodality was first proved by Gasharov \cite{gasharov} whenever 
the rank is at most $2$, and as mentioned by Reiner and Welker 
\cite{reinerwelker} for all graded posets. 

For the relevant definitions concerning the topology behind the 
Charney-Davis Conjecture we refer the reader to 
\cite{charney,reinerwelker,stanleyalg}. 
\begin{conjecture}[Charney-Davis, \cite{charney}]
Let $\Delta$ be a flag simplicial homology $(d-1)$-sphere, where $d$ 
is even. Then the $h$-vector, $h(\Delta,t)$, of $\Delta$ satisfies 
$$
(-1)^{d/2}h(\Delta,-1) \geq 0.
$$
\end{conjecture}
Recall that the $n$th {\em Eulerian polynomial}, $A_n(x)$, is the 
$W$-polynomial of an anti-chain of $n$ elements.  
The Eulerian polynomials can be written as
$$
A_n(x)= \sum_{i=0}^{\lfloor (n-1)/2 \rfloor} a_{n,i}x^i(1+x)^{n-1-2i},
$$
where $a_{n,i}$ is a non-negative integer for all $i$. This was proved 
by Foata and Sch{\"u}tzenberger in \cite{foata} and combinatorially 
by Shapiro, Getu and Woan in \cite{shapiro}. From this expansion 
we see immediately that $A_n(x)$ is symmetric and that the coefficients 
in the standard basis are unimodal. It also follows that 
$(-1)^{(n-1)/2}A_n(-1)\geq 0$.  

We will in Section \ref{signgraded} define a class of labeled 
poset whose members we call sign-graded posets. This class 
includes the class of naturally labeled graded posets. In 
Section \ref{dude} we show that the $W$-polynomial of 
a sign-graded poset $(P,\om)$ of rank $r$ can be expanded, just as the 
Eulerian polynomial, as 
\begin{equation}\label{getu}   
W(P,\om;t) = \sum_{i=0}^{\lfloor (p-r-1)/2 \rfloor}
a_i(P,\om)t^i(1+t)^{p-r-1-2i},
\end{equation}
where $a_i(P,\om)$ are non-negative integers. Hence,  
symmetry and unimodality follow, and $W(P,\om;t)$ 
has the right sign at $-1$. Consequently,   
whenever the associated sphere $\Delta_{eq}(P)$ of a graded poset $P$ is flag 
the Charney-Davis Conjecture holds for $\Delta_{eq}(P)$. We also note 
that all symmetric polynomials with non-positive zeros only, admits 
an expansion such as \eqref{getu}. Hence, that 
$W(P,\om;t)$ has such an expansion can be seen as further 
evidence for the Neggers-Stanley Conjecture.   

In \cite{rsw} the Charney-Davis quantity of a graded naturally labeled 
poset $(P,\om)$ 
of rank $r$ was defined 
to be $(-1)^{(p-1-r)/2}W(P,\om;-1)$. In Section \ref{CDQ} we give 
a combinatorial interpretation of the Charney-Davis quantity as counting 
certain reverse alternating permutations. 
Finally in Section 
\ref{char} we characterize sign-graded posets in terms 
of properties of order polynomials.

\section{Sign-graded posets}\label{signgraded} 
Recall that a poset $P$ is {\em graded} if all maximal chains in $P$ have the 
same length. If $P$ is graded one may associate a {\em rank function} 
$\rho : P \rightarrow \N$ by letting $\rho(x)$ be the length of 
any saturated chain from a minimal element to $x$. The {\em rank} of a 
graded poset $P$ is defined as the length of any maximal chain in $P$. 
In this section we will 
generalize the notion of graded posets to labeled posets. 
  
\begin{figure}\caption{\label{signfig} A sign-graded poset, its two 
labelings and the corresponding rank function.}
\setlength{\unitlength}{10mm}
$$ 
   \vcenter{\xymatrix@R=12pt@C=12pt{
             &&10&&   \\
             &7\e[ru]&6\e[u]&& \\
             &&2\e[u] &9\e[lu]& \\
             && 5\e[u]\e[ru]&4\e[u]& \\ 
             &&& 1\e[u]\e[lu]& 8\\
             &&& 3\e[u]\e[ru]       & 
        }}
\vcenter{\xymatrix@R=20pt@C=20pt{
           &&\n&&   \\
           &\n\e[ru]^1&\n\e[u]_1&& \\
           &&\n\e[u]^1 &\n\e[lu]_{-1}& \\
           && \n\e[u]^{-1}\e[ru]_{1}&\n\e[u]_{1}& \\ 
           &&& \n\e[u]_{1}\e[lu]^1& \n\\
           &&& \n \e[u]^{-1}\e[ru]_{1}       & 
      }}
 \vcenter{\xymatrix@R=12pt@C=12pt{
           &&1&&   \\
           &0\e[ru]&0\e[u]&& \\
           &&-1\e[u] &1 \e[lu]& \\
           && 0\e[u]\e[ru]&0\e[u]& \\ 
           &&& -1\e[u]\e[lu]& 1\\
           &&& 0 \e[u]\e[ru]       & 
      }}
$$
\end{figure}
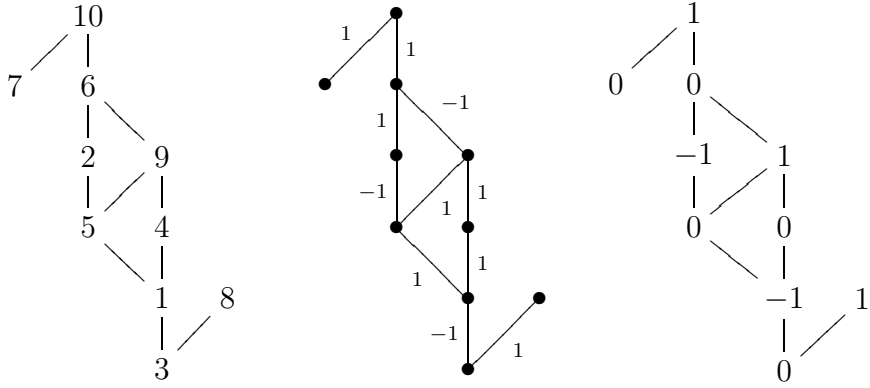

Let $(P,\omega)$ be a labeled poset. An 
element $y$ {\em covers} $x$, written $x \prec y$, if $x<y$ and $x<z<y$ for no 
$z \in P$. Let 
$E=E(P)= \{ (x,y)\in P\times P : x \prec y\}$ be the covering relations 
of $P$. We associate a labeling $\epsilon : E \rightarrow \{-1,1\}$ of the 
covering relations defined by 
$$
\epsilon(x,y)=
\begin{cases} 
\ \ 1  \mbox{ if } \ \ \om(x) < \om(y), \\
-1 \mbox{ if } \ \ \om(x) > \om(y). 
\end{cases}
$$
If two labelings $\om$ and $\lambda$ of $P$ give rise to the same labeling 
of $E(P)$ then it is easy to see that the set of $(P,\om)$-partitions 
and the set of $(P,\lambda)$-partitions are the same.  
In what follows we will often refer to $\epsilon$ as 
the labeling and write $(P,\epsilon)$. 
\begin{definition}
Let $(P,\om)$ be a labeled poset and let $\epsilon$  be the corresponding 
labeling of $E(P)$.  
We say 
that $(P,\om)$ is {\em sign-graded}, and that $P$ is  
$\epsilon$-{\em graded} (and $\om$-{\em graded}) if 
for every maximal chain $x_0 \prec x_1 \prec \cdots \prec x_n$ 
the sum 
$$
\sum_{i=1}^n \epsilon(x_{i-1},x_i)
$$
is the same. The common value of the above sum is called 
the {\em rank} of $(P,\om)$ and is denoted $r(\epsilon)$. 

We say that the poset $P$ is $\epsilon$-{\em consistent} if for 
every $y \in P$ the principal order ideal 
$\Lambda_y = \{x \in P : x \leq y\}$ is $\epsilon_y$-graded, where 
$\epsilon_y$ is $\epsilon$ restricted to $E(\Lambda_y)$.  
The    
{\em rank function} $\rho : P \rightarrow \Z$ of an 
$\epsilon$-consistent poset $P$ is defined by  
$
\rho(x) = r(\epsilon_x). 
$
Hence, an $\epsilon$-consistent poset 
$P$ is $\epsilon$-graded if and only if $\rho$ is constant on 
the set of maximal elements. 
\end{definition}
See Fig. \ref{signfig} for an example of a sign-graded poset. 
Note that if $\epsilon$ is identically equal to $1$, i.e., if $(P,\om)$ is 
naturally labeled, then a sign-graded 
poset with respect to $\epsilon$ is just a graded poset. Note also that 
if $P$ 
is $\epsilon$-graded then $P$ is also $-\epsilon$-graded, where 
$-\epsilon$ is defined by $(-\epsilon)(x,y)= -\epsilon(x,y)$. Up to 
a shift, the order polynomial of a sign-graded labeled poset only depends 
on the underlying poset:  
\begin{theorem}\label{mainman}
Let $P$ be $\epsilon$-graded and $\mu$-graded. Then 
$$
\Omega(P,\epsilon;t-\frac{r(\epsilon)}{2}) = 
\Omega(P,\mu;t-\frac{r(\mu)}{2}).
$$
\end{theorem}
\begin{proof}
Let $\rho_\epsilon$ and $\rho_\mu$ denote the rank functions of
$(P,\epsilon)$ and $(P,\mu)$ respectively, and  
let $\Pa(\epsilon)$ denote the set of 
$(P,\epsilon)$-partitions. Define a function 
$\xi :\Pa(\epsilon) \rightarrow \Q^P$ by 
$
\xi \sigma(x) = \sigma(x) + \Delta(x),
$
where 
$$
\Delta(x)=\frac{r(\epsilon) - \rho_{\epsilon}(x)}{2}
       - \frac{r(\mu) - \rho_{\mu}(x)}{2}.
$$
\begin{table}[h]\label{tab}
\caption{}
{ \begin{tabular}{|l|l|l|l|l|}\hline
$\epsilon(x,y)$ & $\mu(x,y)$ & $\sigma$ & $\Delta$&$\xi\sigma$\\ \hline\hline
$1$ &  $1$  & $\sigma(x) \geq \sigma(y)$& $\Delta(x)=\Delta(y)$& 
$\xi\sigma(x) \geq \xi\sigma(y)$ \\ \hline
$1$ &  $-1$  & $\sigma(x) \geq \sigma(y)$& $\Delta(x)=\Delta(y)+1$& 
$\xi\sigma(x) > \xi\sigma(y)$ \\ \hline
$-1$ &  $1$  & $\sigma(x) > \sigma(y)$& $\Delta(x)=\Delta(y)-1$& 
$\xi\sigma(x) \geq \xi\sigma(y)$ \\ \hline
$-1$ &  $-1$  & $\sigma(x) > \sigma(y)$& $\Delta(x)=\Delta(y)$& 
$\xi\sigma(x) > \xi\sigma(y)$ \\ \hline
\end{tabular}}
\end{table}
The four possible combinations of labelings of a 
covering-relation $(x,y) \in E$ are given in Table 1.   

According to the table $\xi\sigma$ is a $(P,\mu)$-partition provided 
that $\xi\sigma(x)>0$ for all $x \in P$. 
But $\xi\sigma$ is order-reversing so  
it attains 
its minima on maximal elements and if $z$ is a maximal element we 
have $\xi\sigma(z) = \sigma(z)$. Hence  
$\xi : \Pa(\epsilon) \rightarrow \Pa(\mu)$. By symmetry we also have 
a map $\eta : \Pa(\mu) \rightarrow \Pa(\epsilon)$ defined by 
$$
\eta \sigma(x) = \sigma(x) +
        \frac{r(\mu) - \rho_{\mu}(x)}{2}
- \frac{r(\epsilon) - \rho_{\epsilon}(x)}{2}. 
$$
Hence, $\eta= \xi^{-1}$ and $\xi$ is a bijection. 

Since $\sigma$ and $\xi\sigma$ are order-reversing they attain their 
maxima on minimal elements. But if $z$ is a minimal element then  
$\xi\sigma(z) = \sigma(z) + \frac{r(\epsilon)-r(\mu)}{2}$, which gives 
$$
\Omega(P,\mu;n)=\Omega(P,\epsilon;n + \frac{r(\mu)-r(\epsilon)}{2}),
$$
for all non-negative integers $n$ and the theorem follows.
\end{proof}

\begin{theorem}\label{ordersym}
Let $P$ be $\epsilon$-graded. Then 
$$
\Omega(P,\epsilon;t) = (-1)^{p}\Omega(P,\epsilon;-t-r(\epsilon)).
$$ 
\end{theorem}
\begin{proof}
We have the following reciprocity for order polynomials, see 
\cite{stanleythesis}:
\begin{equation}\label{orderrec}
\Omega(P,-\epsilon;t) = (-1)^{p}\Omega(P,\epsilon;-t).
\end{equation} 
Note that $r(-\epsilon)=-r(\epsilon)$, so by Theorem \ref{mainman} we 
have:
$$
\Omega(P,-\epsilon;t)=\Omega(P,\epsilon,t-r(\epsilon)),
$$
which, combined with \eqref{orderrec}, gives the desired result.
\end{proof}

\begin{corollary}\label{symmetri}
Let $P$ be an $\epsilon$-graded poset. Then $W(P, \epsilon, t)$ is 
symmetric with center of symmetry $(p-r(\epsilon)-1)/2$. 
If $P$ is also $\mu$-graded then 
$$
W(P,\mu;t)= t^{(r(\epsilon)-r(\mu))/2}W(P,\epsilon;t).
$$  
\end{corollary}
\begin{proof}
It is known, see \cite{stanleythesis}, that if $W(P,\epsilon;t)= 
\sum_{i \geq 0} w_i(P,\epsilon)t^i$ then 
$\Omega(P,\epsilon;t) = \sum_{i \geq 0}w_i(P,\epsilon)\binom {t+p-1-i}p$. 
Let $r=r(\epsilon)$.  
Theorem \ref{ordersym} gives:
\begin{eqnarray*}
\Omega(P,\epsilon;t) &=& 
\sum_{i \geq 0} w_i(P,\epsilon) (-1)^p\binom{-t-r+p-1-i}p \\
                &=& \sum_{i \geq 0} w_i(P,\epsilon) \binom{t+r+i}p \\
                &=& \sum_{i \geq 0} w_{p-r-1-i}(P,\epsilon) \binom{t+p-1-i}p,
\end{eqnarray*} 
so $w_i(P,\epsilon) = w_{p-r-1-i}(P,\epsilon)$ for all $i$, and the symmetry 
follows. The relationship between the $W$-polynomials of $\epsilon$ and 
$\mu$ follows from 
Theorem \ref{mainman} and the expansion of order-polynomials in the 
basis $\binom {t+p-1-i}p$.  
\end{proof}
We say that a poset $P$ is {\em parity graded} if the size of  
all maximal chains in $P$ have the same parity. Also, a poset is $P$ is 
{\em parity consistent} if for all $x \in P$ the order ideal 
$\Lambda_x$ is parity graded. These classes of posets were 
studied in \cite{stanleybalance} in a different context.   
The following theorem tells us that the class of sign-graded posets 
is considerably greater than the class of graded posets.
\begin{theorem}\label{evenodd}
Let $P$ be a poset. Then 
\begin{itemize}
\item there exists a labeling 
$\epsilon : E \rightarrow \{-1,1\}$ such that $P$ is 
$\epsilon$-consistent if and only if $P$ is parity consistent,
\item there exists a labeling 
$\epsilon : E \rightarrow \{-1,1\}$ such that $P$ is  
$\epsilon$-graded if and only if $P$ is parity graded.
\end{itemize} 
Moreover, the labeling $\epsilon$ can be chosen so 
that the corresponding rank function has values in $\{0,1\}$.
\end{theorem}
\begin{proof}
It suffices to prove the equivalence regarding parity graded posets.  
It is clear that if $P$ is $\epsilon$-graded then $P$ is 
parity graded. 

Let $P$ be parity graded. Then, for any $x \in P$, 
all saturated from a minimal 
element to $x$ has the same length modulo $2$.  
Hence, we may define a labeling $\epsilon: P \rightarrow \{-1,1\}$ 
by $\epsilon(x,y)= (-1)^{\ell(x)}$, where $\ell(x)$ is the length of 
any saturated chain starting at a minimal element and ending at $x$. 
It follows that $P$ is $\epsilon$-graded and that its  
rank function has values in $\{0,1\}$.
\end{proof} 
We say that $\om : P \rightarrow \{1,2,\ldots ,p\}$ is 
{\em canonical} if  $(P,\om)$ has a rank-function $\rho$ with values in 
$\{0,1\}$, and $\rho(x)<\rho(y)$ implies $\om(x) < \om(y)$. By 
Theorem \ref{evenodd} we know that $P$ admits a canonical labeling 
if $P$ is $\epsilon$-consistent for some $\epsilon$.
\section{The Jordan-H\"older set of an $\epsilon$-consistent poset}
Let $P$ be $\omega$-consistent.  
We may assume that $\omega(x)<\omega(y)$ 
whenever $\rho(x)<\rho(y)$. Suppose that 
$x,y \in P$ are incomparable and that $\rho(y)=\rho(x)+1$. Then  
the Jordan-H\"older set of $(P,\omega)$ can be partitioned into 
two sets: One where in all permutations $\om(x)$ comes before $\om(y)$ and 
one where $\om(y)$ always comes before $\om(x)$. This means that 
$\J(P,\omega)$ is the disjoint union
\begin{equation}\label{split}
\J(P,\omega)= \J(P',\omega) \sqcup \J(P'',\omega),
\end{equation}
where $P'$ is the transitive closure of $E\cup \{x \prec y\}$, 
and $P''$ is the transitive closure 
of $E \cup \{y \prec x\}$.
\begin{lemma}\label{primes}
With definitions as above $P'$ and $P''$ are 
$\omega$-consistent with the same rank-function as $(P,\omega)$.
\end{lemma}
\begin{proof}
Let 
$
c: z_0 \prec z_1 \prec \cdots \prec z_k=z
$
be a saturated chain in $P''$, where $z_0$ is a minimal element of $P''$. Of 
course $z_0$ is also a minimal element of $P$.  
We have to prove that 
$$
\rho(z) = \sum_{i=0}^{k-1}\epsilon''(z_i,z_{i+1}),
$$
where $\epsilon''$ is the labeling of $E(P'')$ and $\rho$ is 
the rank-function of $(P,\omega)$. 

All covering relations in $P''$, except $y \prec x$, are also covering 
relations in $P$. If $y$ and $x$ do not 
appear in $c$, then $c$ is a saturated chain in $P$ and there is  nothing 
to prove. Otherwise 
$$
c: y_0 \prec \cdots \prec y_i=y \prec x=x_{i+1} \prec x_{i+2} 
\prec \cdots \prec x_k=z.
$$
Note that if $s_0 \prec s_1 \prec \cdots \prec s_\ell$ is any saturated 
chain in $P$ then 
$\sum_{i=0}^{\ell-1}\epsilon(s_i,s_{i+1}) = \rho(s_\ell)-\rho(s_0)$. 
Since $y_0 \prec \cdots \prec y_i=y$ and $x=x_{i+1} \prec x_{i+2} 
\prec \cdots \prec x_k=z$ are saturated chains in $P$ we have 
\begin{eqnarray*}
\sum_{i=0}^{k-1}\epsilon''(z_i,z_{i+1})&=&
\rho(y) + \epsilon''(y,x)+ \rho(z)-\rho(x)\\
                                   &=& \rho(y)-1-\rho(x)+\rho(z)\\
&=&\rho(z),
\end{eqnarray*}
as was to be proved. The statement for $(P',\omega)$ follows similarly.   
\end{proof}
We say that a $\omega$-consistent poset $P$ is {\em saturated}  
if for all $x, y \in P$ we have that $x$ and $y$ are comparable whenever 
$|\rho(y)-\rho(x)|=1$. Let $P$ and $Q$ be posets on the same set. Then 
$Q$ {\em extends} $P$ if $x <_Q y$ whenever $x <_P y$. 
\begin{theorem}\label{saturated}
Let $P$ be a $\omega$-consistent poset. Then 
the Jordan-H\"older set of $(P,\omega)$ is uniquely decomposed as the 
disjoint union
$$
\J(P,\omega)= \bigsqcup_{Q} \J(Q,\omega),
$$
where the union is over all saturated $\omega$-consistent posets 
$Q$   
that extend $P$ and have the same rank-function as $(P, \omega)$. 
\end{theorem}
\begin{proof}
That the union exhausts $\J(P,\omega)$ follows from \eqref{split} and 
Lemma \ref{primes}. Let  
$Q_1$ and $Q_2$ be two different saturated $\omega$-consistent 
posets that extend $P$ and have the same rank-function as 
$(P,\om)$. We may assume that $Q_2$ does not extend $Q_1$. Then 
there exists a covering relation $x \prec y$ in $Q_1$ such that 
$x \nless y$ in $Q_2$.  Since 
$|\rho(x) - \rho(y)|=1$ we must have $y < x$ in $Q_2$. Thus 
$\om(x)$ precedes $\om(y)$ in any permutation in $\J(Q_1,\om)$, 
and  $\om(y)$ precedes $\om(x)$ in any permutation in $\J(Q_2,\om)$. Hence, 
the union is disjoint and unique.   
\end{proof}
We need two operations on labeled posets:
Let $(P,\epsilon)$ and 
$(Q,\mu)$ be two labeled posets. 
The {\em ordinal sum}, $P \oplus Q$, of  $P$ and $Q$ 
is the poset with 
the disjoint union of $P$ and $Q$ as underlying set and with partial order 
defined by $x \leq y$ if $x \leq_P y$ or $x \leq_Q y$, or 
$x \in P, y \in Q$. Define two labelings of $E(P \oplus Q)$ by
\begin{eqnarray*}
(\epsilon \oplus_1 \mu)(x,y)&=& \epsilon(x,y) \mbox{ if } (x,y) \in E(P),\\ 
(\epsilon \oplus_1 \mu)(x,y)&=& \mu(x,y) \mbox{ if } (x,y) \in E(Q) 
\mbox{ and } \\ 
(\epsilon \oplus_1 \mu)(x,y)&=& 1  \mbox{ otherwise.}\\   
(\epsilon \oplus_{-1} \mu)(x,y)&=&\epsilon(x,y) \mbox{ if } (x,y) \in E(P),\\ 
(\epsilon \oplus_{-1} \mu)(x,y)&=&\mu(x,y) \mbox{ if } (x,y) \in E(Q)  
\mbox{ and } \\ 
(\epsilon \oplus_{-1} \mu)(x,y)&=&-1 \mbox{ otherwise.}
\end{eqnarray*}
With a slight abuse of notation we write $P\oplus_{\pm 1} Q$ when the 
labelings of $P$ and $Q$ are understood from the context.  
Note that ordinal sums are associative, i.e., 
$(P\oplus_{\pm 1} Q)\oplus_{\pm 1} R= P\oplus_{\pm 1} (Q\oplus_{\pm 1} R)$, 
and preserve the  
property of being sign-graded. 
The following result is easily obtained by combinatorial reasoning, see 
\cite{brentithesis, wagnernc}:
\begin{proposition}\label{ordinalsums}
Let $(P,\om)$ and $(Q,\nu)$ be two labeled posets. Then 
$$
W(P\oplus Q, \om \oplus_1 \nu;t)= W(P,\om;t)W(Q,\nu;t)
$$
and
$$ 
W(P\oplus Q, \om \oplus_{-1} \nu;t)= tW(P,\om;t)W(Q,\nu;t).
$$
\end{proposition}

\begin{proposition}\label{gammeldansk}
Suppose that $(P,\omega)$ is a saturated canonically labeled 
$\omega$-consistent poset. Then 
$(P,\om)$ is the direct sum 
$$
(P,\om) = A_0 \oplus_1 A_1 \oplus_{-1} A_2 \oplus_1 A_3 \oplus_{-1}\cdots 
          \oplus_{\pm 1}A_k,
$$
where the $A_i$s are anti-chains.  
\end{proposition}  
\begin{proof}
Let $\pi \in \J(P,\om)$. Then  we may write $\pi$ as 
$\pi = w_0w_1\cdots w_k$ where 
the $w_i$s are maximal words with respect to the property: 
If $a$ and $b$ are letters of $w_i$ then 
$\rho(\om^{-1}(a))= \rho(\om^{-1}(b))$. Hence $\pi \in J(Q,\om)$ where 
$$
(Q,\om) = A_0 \oplus_1 A_1 \oplus_{-1} A_2 \oplus_1 A_3 \oplus_{-1}
          \cdots \oplus_{\pm 1}A_k,
$$
and $A_i$ is the anti-chain consisting of the elements $\om^{-1}(a)$, 
where $a$ is a letter of $w_i$ ($A_i$ is an anti-chain, since 
if $x<y$ where $x,y \in A_i$ there would be a letter in $\pi$ between 
$\om(x)$ and $\om(y)$ whose rank was different than that of $x,y$).  
Now, $(Q,\om)$ is 
saturated so $P=Q$.   
\end{proof}
Note that the argument in the above proof also can be used to 
give a simple proof of Theorem \ref{saturated} when $\om$ is 
canonical. 

\section{The $W$-polynomial of a sign-graded poset}\label{dude}
The space $S^d$  of symmetric polynomials in $\R[t]$ with center 
of symmetry $d/2$ has a basis 
$$
B_d=\{t^i(1+t)^{d-2i}\}_{i=0}^{\lfloor d/2 \rfloor}.
$$ 
If $h \in S^d$ has non-negative coefficients in this basis it follows 
immediately that the coefficients of $h$ in the standard basis are unimodal. 
Let $S_+^d$ be the non-negative span of $B_d$. 
Thus $S_+^d$ is a cone. Another property of $S_+^d$ is that 
if $h \in S_+^d$ then it has the correct sign at $-1$ i.e., 
$$
(-1)^{d/2}h(-1) \geq 0.
$$ 
\begin{lemma}\label{ssp}
Let $c,d \in \N$. Then 
\begin{eqnarray*}
S^cS^d &\subset& S^{c+d} \\
S_+^cS_+^d &\subset& S_+^{c+d}.
\end{eqnarray*}
Suppose further that  $h \in S^d$ has positive leading coefficient 
and that all zeros of $h$ are real and non-positive. Then $h \in S_+^d$. 
\end{lemma}
\begin{proof}
The inclusions are obvious. Since $t \in S_+^2$ and $(1+t) \in S_+^1$ we may 
assume that
none of them divides $h$. But then we may collect the zeros of $h$ in 
pairs $\{ \theta,\theta^{-1}\}$. Let $A_\theta = -\theta -\theta^{-1}$. Then 
$$
h= C\prod_{\theta<-1}(t^2 + A_\theta t +1),
$$
where $C>0$. Since $A_\theta >2$ we have 
$$
t^2 + A_\theta t +1= (t+1)^2 + (A_\theta -2)t \in S_+^2,
$$
and the lemma follows.
\end{proof} 

We can now prove our main theorem.
\begin{theorem}\label{mainthm}
Suppose that $(P,\om)$ is a sign-graded poset of rank $r$. 
Then 
$W(P,\om;t) \in S_+^{p-r-1}$. 
\end{theorem}
\begin{proof}
By Corollary \ref{symmetri} and Lemma \ref{evenodd} we may assume 
that $(P,\om)$ is canonically labeled. If $Q$ extends $P$ then 
the maximal elements of $Q$ are also maximal elements of $P$.     
By Theorem \ref{saturated} we know that 
$$
W(P,\om;t) = \sum_{Q} W(Q,\om;t),
$$
where $(Q,\om)$ is saturated and sign-graded with the 
same rank function and rank as $(P,\omega)$. 
The $W$-polynomials of anti-chains are the Eulerian polynomials, which 
have real non-negative zeros only.  
By Propositions \ref{ordinalsums} and 
\ref{gammeldansk} the polynomial 
$W(Q,\om;t)$ has only real non-positive zeros so by 
Lemma \ref{ssp} and Corollary \ref{symmetri} we have 
$W(Q,\om;t) \in S_+^{p-r-1}$. The 
theorem now follows since $S_+^{p-r-1}$ is a cone.
\end{proof}
\begin{corollary}
Let $(P,\om)$ be sign-graded of rank $r$. Then 
$W(P,\om;t)$ is symmetric and its coefficients are unimodal. Moreover, 
$W(P,\om;t)$ has the correct sign at $-1$, i.e.,
$$
(-1)^{(p-1-r)/2}W(P,\om;-1) \geq 0.
$$
\end{corollary}
\begin{corollary}
Let $P$ be a graded poset. Suppose that 
$\Delta_{eq}(P)$ is flag. Then the Charney-Davis Conjecture 
holds for $\Delta_{eq}(P)$.
\end{corollary}

\begin{theorem}
Suppose that $P$ is an  
$\om$-consistent poset and that $|\rho(x)-\rho(y)|\leq 1$ for all 
 maximal 
elements $x,y \in P$. Then $W(P,\om;t)$ has unimodal 
coefficients.
\end{theorem}
\begin{proof}
Suppose that the ranks of maximal elements are contained in $\{r,r+1\}$. 
If $Q$ is any saturated poset that extends $P$ and has the same rank 
function as $(P,\om)$ then $Q$ is $\om$-graded of rank $r$ or $r+1$.
By Theorems \ref{saturated} and \ref{mainthm} we know that 
$$
W(P,\om;t)=\sum_{Q}W(Q,\om;t),
$$
where $W(Q,\om;t)$ is symmetric and unimodal with center of symmetry 
at $(p-1-r)/2$ or $(p-2-r)/2$. The 
sum of such polynomials is again unimodal. 
\end{proof}
\section{The Charney-Davis quantity}\label{CDQ}
In \cite{rsw} Reiner, Stanton and Welker defined the {\em Charney-Davis 
quantity} of a graded naturally labeled poset $(P,\om)$ of rank $r$ to 
be 
$$
CD(P,\om)=(-1)^{(p-1-r)/2}W(P,\om;-1).
$$
We define it in the exact same way for sign-graded posets. Since 
the particular labeling does not matter we write $CD(P)$. Let 
$\pi=\pi_1\pi_2 \cdots \pi_n$ be any permutation. We say that 
$\pi$ is {\em alternating} if $\pi_1>\pi_2<\pi_3> \cdots$ and   
{\em reverse alternating} if $\pi_1<\pi_2>\pi_3< \cdots$. 
Let $(P,\om)$ be a canonically labeled sign-graded poset. If 
$\pi \in \J(P,\om)$ then we may write $\pi$ as 
$\pi = w_0w_1\cdots w_k$ where $w_i$ are maximal words with respect 
to the property: If $a$ and $b$ are letters of $w_i$ then 
$\rho(\om^{-1}(a))=\rho(\om^{-1}(b))$. The words $w_i$ are called 
the {\em components} of $\pi$. The  following theorem  is well known, 
see for example \cite{shapiro}, and gives the Charney-Davis quantity 
of an anti-chain.
\begin{proposition}\label{symalt}
Let $n \geq 0$ be an integer. Then 
$
(-1)^{(n-1)/2}A_n(-1)
$
is equal to $0$ if $n$ is even and equal to the number of (reverse) 
alternating permutations of the set $\{1,2,\ldots,n\}$ if $n$ is odd.
\end{proposition}  

\begin{theorem}\label{cdq}
Let $(P,\om)$ be a canonically labeled sign-graded poset. Then the 
Charney-Davis quantity, $CD(P)$, is equal to the number of reverse 
alternating permutations in $\J(P,\om)$ such that all 
components have an odd numbers of letters. 
\end{theorem}
\begin{proof}
It suffices to prove the theorem when $(P,\om)$ is saturated. By Proposition 
\ref{gammeldansk} we know that 
$$
(P,\om) = A_0 \oplus_1 A_1 \oplus_{-1} A_2 \oplus_1 A_3 \oplus_{-1}\cdots 
\oplus_{\pm 1}A_k, 
$$
where the $A_i$s are anti-chains. Thus  
$CD(P)= CD(A_0)CD(A_1)\cdots CD(A_k)$. 
Let $\pi =w_0w_1 \cdots w_k \in \J(P,\om)$ where $w_i$ is a 
permutation of $\om(A_i)$. Then $\pi$ is a reverse 
alternating such that all 
components have an odd numbers of letters if and only if, for all $i$,  
$w_i$ is reverse alternating 
if $i$ is even and alternating if $i$ is odd. Hence, by Proposition 
\ref{symalt}, the 
number of such permutations is indeed $CD(A_0)CD(A_1)\cdots CD(A_k)$.
\end{proof} 
If $h(t)$ is any polynomial with integer coefficients and 
$h(t) \in S^d$, it follows that $h(t)$ has integer 
coefficients in the basis $t^i(1+t)^{d-2i}$. Thus we know 
that if $(P,\om)$ is sign-graded of rank $r$, then 
$$
W(P,\om;t) = \sum_{i=0}^{\lfloor (p-r-1)/2 \rfloor}
              a_i(P,\om)t^i(1+t)^{p-r-1-2i},
$$
where $a_i(P,\om)$ are non-negative integers. By Theorem 
\ref{cdq} we have a combinatorial interpretation of 
the $a_{(p-r-1)/2}(P,\om)$. A similar but more complicated interpretation of  
$a_i(P,\om)$, $i=0,1,\ldots,\lfloor (p-r-1)/2 \rfloor$ can be deduced from 
Proposition \ref{gammeldansk} and the work in \cite{shapiro}. We omit 
this.   
\section{The right mode}
Let $f(x)=a_0 + a_1x + \cdots + a_dx^d$ be a polyomial with real 
coefficients. The {\em mode}, $\mode(f)$, of $f$ is the average value 
of the indices $i$ such that $a_i=\max\{a_j\}_{j=0}^d$. One can easily compute 
the mode of a polynomial with real non-positive zeros only:
\begin{theorem}\cite{darroch}\label{rightmode}
Let $f$ be a polynomial with real non-positive zeros only and with 
positive leading coefficient. Then 
$$
\Big| \frac {f'(1)}{f(1)} - \mode(f) \Big|<1.
$$
\end{theorem} 
It is known, see 
\cite{brentithesis,stanleythesis,wagnernc}, that  
$$
W(P,\om;x) = \sum_{i=1}^pe_i(P,\om)x^{i-1}(1-x)^{p-i}, 
$$
where $e_i(P,\om)$ is the number of surjective $(P,\om)$-partitions 
$\sigma : P \rightarrow \{1,2,\ldots,i\}$.
A simple calculation gives
\begin{equation}\label{mode}
\frac {W'(P,\om;1)} {W(P,\om;1)} = p-1 - \frac {e_{p-1}(P,\om)}{e_p(P,\om)}.
\end{equation} 
If $P$ is $\om$-graded of rank $r$ we know by Theorem \ref{mainthm} that 
$\mode(W(P,\om;x))= (p-r-1)/2$. The Neggers-Stanley conjecture, 
Theorem \ref{rightmode} and 
\eqref{mode} suggest  
that $2e_{p-1}(P,\om)=(p+r-1)e_p(P,\om)$. Stanley \cite{stanleythesis} 
proved this for graded posets and it generalizes to sign-graded posets:
\begin{proposition}
Let $P$ be $\om$-graded of rank $r$. Then 
$$
2e_{p-1}(P,\om)=(p+r-1)e_p(P,\om).
$$
\end{proposition}
\begin{proof}
The identity follows when expanding $\Omega(P,\om,t)$ in powers of $t$ using 
Theorem \ref{ordersym}. See \cite[Corollary 19.4]{stanleythesis} for 
details.
\end{proof}

\section{A characterization of sign-graded posets}\label{char}
Here we give a characterization of sign-graded posets along 
the lines of the characterization of graded posets given by 
Stanley in \cite{stanleythesis}.
Let $(P,\epsilon)$ be a labeled poset. Define a function 
$\delta=\delta_\epsilon : P \rightarrow \Z$ by 
$$
\delta(x) = \max \{ \sum_{i=1}^{\ell} \epsilon(x_{i-1},x_i) \},
$$
where $x=x_0 \prec x_1 \prec \cdots \prec x_{\ell}$ is any saturated chain 
starting at $x$ and ending at a maximal element $x_\ell$. Define 
a map $\Phi=\Phi_\epsilon : \Pa(\epsilon) \rightarrow \Z^P$ by 
$$
\Phi\sigma = \sigma + \delta. 
$$
We have  
\begin{equation}\label{deltaineq}
\delta(x) \geq 
\delta(y) +\epsilon(x,y).
\end{equation}
This means that  
$\Phi\sigma(x) > \Phi\sigma(y)$ if $\epsilon(x,y)=1$ and  
$\Phi\sigma(x) \geq \Phi\sigma(y)$ if $\epsilon(x,y)=-1$. Thus $\Phi\sigma$ 
is a $(P,-\epsilon)$-partition provided that $\Phi \sigma(x)>0$ for 
all $x \in P$. But $\Phi\sigma$ is order reversing so it attains 
its minimum at maximal elements and for maximal elements, $z$, we 
have $\Phi\sigma(z) = \sigma(z)$. This shows that 
$\Phi : \Pa(\epsilon) \rightarrow \Pa(-\epsilon)$ is an injection. 

The {\em dual}, $(P^*, \epsilon^*)$, of a labeled poset 
$(P,\epsilon)$ is defined by $x <_{P^*} y$ if and only if $y <_{P^*} x$, 
with labeling defined by $\epsilon^*(y,x)= -\epsilon(x,y)$. We say that 
$P$ is {\em dual} $\epsilon$-{\em consistent} if 
$P^*$ is $\epsilon^*$-consistent.  
\begin{proposition}\label{delta}
Let $(P,\epsilon)$ be labeled poset. Then 
$\Phi_\epsilon : \Pa(\epsilon) \rightarrow \Pa(-\epsilon)$  is a 
bijection if and only if $P$ is dual $\epsilon$-consistent. 
\end{proposition}
\begin{proof}
If $P$ is dual $\epsilon$-consistent then $P$ is also dual 
$-\epsilon$-consistent and $\delta_{-\epsilon}(x) = -\delta_\epsilon(x)$ for 
all $x \in P$. Thus the if part follows since the inverse of 
$\Phi_\epsilon$ is $\Phi_{-\epsilon}$. 

For the only if direction note that $P$ is dual $\epsilon$-consistent 
if and only if for all $(x,y) \in E$ we have 
$$
\delta(x) = 
\delta(y) + \epsilon(x,y)
$$
Hence, if $P$ is not dual $\epsilon$-consistent then 
by \eqref{deltaineq}, there is a covering relation $(x_0,y_0) \in E$ such 
that either $\epsilon(x_0,y_0)=1$ and $\delta(x_0) \geq \delta(y_0)+2$ or 
$\epsilon(x_0,y_0)=-1$ and $\delta(x_0) \geq \delta(y_0)$. 

Suppose that $\epsilon(x_0,y_0)=1$. It is clear that there is a 
$\sigma \in \Pa(-\epsilon)$ such that $\sigma(x_0)=\sigma(y_0)+1$. But then 
$$
\sigma(x_0) - \delta(x_0) \leq \sigma(y_0)-\delta(y_0)-1, 
$$
so $\sigma - \delta \notin \Pa(\epsilon)$. 

Similarly, if $\epsilon(x_0,y_0)=-1$ then we can find a partition 
$\sigma \in \Pa(-\epsilon)$ with $\sigma(x_0)=\sigma(y_0)$, and then  
$$
\sigma(x_0) - \delta(x_0) \leq \sigma(y_0)-\delta(y_0), 
$$
so $\sigma - \delta \notin \Pa(\epsilon)$.  
\end{proof}
Let $(P,\epsilon)$ be a labeled poset. 
Define $r(\epsilon)$ by 
$$
r(\epsilon) = \max  \{ \sum_{i=1}^{\ell} \epsilon(x_{i-1},x_i) : 
x_0 \prec x_1 \prec \cdots \prec x_{\ell} \mbox{ is maximal}\}.  
$$
We then have: 
\begin{eqnarray*}
\max \{ \Phi\sigma(x) : x \in P\} &=& 
\max \{ \sigma(x) + \delta_\epsilon(x): x \mbox{ is minimal}\} \\
&\leq& \max \{ \sigma(x) : x \in P\} + r(\epsilon). 
\end{eqnarray*}
So if we let $\Pa_n(\epsilon)$ be the $(P,\epsilon)$-partitions with 
largest part at most $n$ we have that 
$\Phi_\epsilon : \Pa_n(\epsilon) \rightarrow \Pa_{n+r(\epsilon)}(-\epsilon)$
is an injection. A labeling $\epsilon$ of $P$ is said to satisfy the 
$\lambda$-{\em chain condition} if for every $x \in P$ there is 
a maximal chain $c : x_0 \prec x_1 \prec \cdots \prec x_\ell$ containing 
$x$ such that $\sum_{i=1}^\ell \epsilon(x_{i-1},x_i)= r(\epsilon)$. 
\begin{lemma}\label{lambda}
Suppose that $n$ is a non-negative integer such that 
$\Omega(P,\epsilon;n) \neq 0$. If 
$$
\Omega(P,-\epsilon;n+r(\epsilon)) = \Omega(P, \epsilon;n)
$$ 
then $\epsilon$ satisfies the $\lambda$-chain condition. 
\end{lemma}
\begin{proof}
Define $\delta^* : P \rightarrow \Z$ by 
$$
\delta^*(x) = \max\{ \sum_{i=1}^\ell\epsilon(x_{i-1},x_i): 
x_0 \prec x_1 \prec \cdots \prec x_{\ell}=x  \},
$$
where the maximum is taken over all maximal chains starting at a 
minimal element and ending at $x$. Then 
\begin{equation}\label{star}
\delta(x) + \delta^*(x) \leq r(\epsilon)
\end{equation}
for all $x$, and $\epsilon$ satisfies the $\lambda$-chain condition if 
and only if we have equality in \eqref{star} for all $x \in P$. It is 
easy to see that the map 
$\Phi^* : \Pa_n(\epsilon) \rightarrow \Pa_{n+r(\epsilon)}(-\epsilon)$ defined 
by 
$$
\Phi^*\sigma(x) = \sigma(x) +r(\epsilon) - \delta^*(x), 
$$
is well-defined and is an injection. By \eqref{star} we have 
$\Phi\sigma(x) \leq \Phi^*\sigma(x)$ for all 
$\sigma$ and all $x \in P$, with equality if and only if $x$ is in a 
maximal chain of maximal weight. This means that in order for 
$\Phi : \Pa_n(\epsilon) \rightarrow \Pa_{n+r(\epsilon)}(-\epsilon)$ to 
be a bijection it is necessary for $\epsilon$ to satisfy the 
$\lambda$-chain condition.     
\end{proof}  
\begin{theorem}\label{omm}
Let $\epsilon$ be a labeling of $P$. Then 
$$
\Omega(P,\epsilon;t) = (-1)^{p}\Omega(P,\epsilon;-t-r(\epsilon))
$$
if and only if $P$ is $\epsilon$-graded of rank $r(\epsilon)$. 
\end{theorem}
\begin{proof}
The ''if'' part is Theorem \ref{ordersym}, so suppose that the 
equality of the theorem holds. By reciprocity we have 
$$
(-1)^{p}\Omega(P,\epsilon;-t-r(\epsilon))=\Omega(P,-\epsilon;t+r(\epsilon)),
$$
and since $\Phi_\epsilon : \Pa_n(\epsilon) \rightarrow 
\Pa_{n+r(\epsilon)}(-\epsilon)$ is an injection it is also a bijection. 
By Proposition \ref{delta} we have that $P$ is dual 
$\epsilon$-consistent and by Lemma \ref{lambda}, 
we have that all minimal elements are members of maximal chains of 
maximal weight. In other words $P$ is $\epsilon$-graded. 
\end{proof}
It should be noted that it is not necessary for $P$ to be 
$\epsilon$-graded in order for $W(P,\epsilon;t)$ to be symmetric. 
For example,  
if $(P, \epsilon)$ is any labeled poset then the $W$-polynomial 
of the disjoint union of $(P,\epsilon)$ and $(P,-\epsilon)$ is easily 
seen to be symmetric. However, we have the following:
\begin{corollary}
Suppose that 
$$
\Omega(P,\epsilon;t)=\Omega(P,-\epsilon;t+s), 
$$ 
for some $s \in \Z$. Then $-r(-\epsilon) \leq s \leq r(\epsilon)$, 
with equality if and only if $P$ is $\epsilon$-graded. 
\end{corollary}
\begin{proof}
We have an injection 
$\Phi_\epsilon : \Pa_n(\epsilon) \rightarrow \Pa_{n+r(\epsilon)}(-\epsilon)$. 
This means that $s \leq  r(\epsilon)$. The lower bound follows from 
the injection $\Phi_{-\epsilon}$, and the statement of equality follows 
from Theorem \ref{omm}. 
\end{proof}

\bibliography{signgraded}

\begin{thebibliography}{10}

\bibitem{branden1}
P.~Br\"and\'en.
\newblock On operators on polynomials preserving real-rootedness and the
  {N}eggers-{S}tanley conjecture.
\newblock {\em J.Algebraic Comb.}, to appear.

\bibitem{brentithesis}
F.~Brenti.
\newblock Unimodal, log-concave and {P}\'olya frequency sequences in
  combinatorics.
\newblock {\em Mem. Amer. Math. Soc.}, 81(413):viii+106, 1989.

\bibitem{charney}
R.~Charney and M.~Davis.
\newblock The {E}uler characteristic of a nonpositively curved, piecewise
  {E}uclidean manifold.
\newblock {\em Pacific J. Math.}, 171(1):117--137, 1995.

\bibitem{darroch}
J.~N. Darroch.
\newblock On the distribution of the number of successes in independent trials.
\newblock {\em Ann. Math. Statist.}, 35:1317--1321, 1964.

\bibitem{foata}
D.~Foata and M.~Sch{\"u}tzenberger.
\newblock {\em Th\'eorie g\'eom\'etrique des polyn\^omes eul\'eriens}.
\newblock Lecture Notes in Mathematics, Vol. 138. Springer-Verlag, Berlin,
  1970.

\bibitem{gasharov}
V.~Gasharov.
\newblock On the {N}eggers-{S}tanley conjecture and the {E}ulerian polynomials.
\newblock {\em J. Combin. Theory Ser. A}, 82(2):134--146, 1998.

\bibitem{neggers}
J.~Neggers.
\newblock Representations of finite partially ordered sets.
\newblock {\em J. Combin. Inform. System Sci.}, 3(3):113--133, 1978.

\bibitem{rsw}
V.~Reiner, D.~Stanton, and V.~Welker.
\newblock The {C}harney-{D}avis quantity for certain graded posets.
\newblock {\em S\'em. Lothar. Combin.}, 50:Art. B50c, 13 pp. (electronic),
  2003.

\bibitem{reinerwelker}
V.~Reiner and V.~Welker.
\newblock On the {C}harney-{D}avis and the {N}eggers-{S}tanley conjectures.
\newblock {\em http://www.math.umn.edu/\~{}reiner/Papers/papers.html}, 2002.

\bibitem{shapiro}
L.~W. Shapiro, W.~J. Woan, and S.~Getu.
\newblock Runs, slides and moments.
\newblock {\em SIAM J. Algebraic Discrete Methods}, 4(4):459--466, 1983.

\bibitem{stanleybalance}
R.~P. Stanley.
\newblock {Some remarks on sign-balanced and maj-balanced posets}.
\newblock {\em arXiv:math.CO/0211113}.

\bibitem{stanleythesis}
R.~P. Stanley.
\newblock {\em Ordered structures and partitions}.
\newblock American Mathematical Society, Providence, R.I., 1972.
\newblock Memoirs of the American Mathematical Society, No. 119.

\bibitem{stanleyg}
R.~P. Stanley.
\newblock The number of faces of simplicial polytopes and spheres.
\newblock In {\em Discrete geometry and convexity (New York, 1982)}, volume 440
  of {\em Ann. New York Acad. Sci.}, pages 212--223. New York Acad. Sci., New
  York, 1985.

\bibitem{stanleyalg}
R.~P. Stanley.
\newblock {\em Combinatorics and commutative algebra}, volume~41 of {\em
  Progress in Mathematics}.
\newblock Birkh\"auser Boston Inc., Boston, MA, second edition, 1996.

\bibitem{stanley1}
R.~P. Stanley.
\newblock {\em Enumerative combinatorics. {V}ol. 1}, volume~49 of {\em
  Cambridge Studies in Advanced Mathematics}.
\newblock Cambridge University Press, Cambridge, 1997.

\bibitem{wagnernc}
D.~G. Wagner.
\newblock Enumeration of functions from posets to chains.
\newblock {\em European J. Combin.}, 13(4):313--324, 1992.

\end{thebibliography}

\end{document}